\documentclass[12pt,reqno]{amsart}
\usepackage{latexsym,amsmath,amsfonts,amssymb,amsthm}
\theoremstyle{plain}
\newtheorem{Thm}{Theorem}
\newtheorem{Lem}{Lemma}
\newtheorem{Cor}{Corollary}
\newtheorem{Conj}{Conjecture}
\theoremstyle{definition}

\theoremstyle{remark}
\newtheorem*{Rem}{Remark}
\def\Z{\mathbb Z}

\def\C{\mathbb C}
\def\P{\mathcal P}
\def\pmod #1{\ ({\rm mod}\ #1)}

\def\supp{{\rm supp}}
\def\floor #1{\lfloor{#1}\rfloor}
\def\ceil #1{\left\lceil{#1}\right\rceil}
\def\1{\mathbf 1}

\begin{document}
\title{A Schur-type addition theorem for primes}
\author{Hongze Li}
\email{lihz@sjtu.edu.cn}
\author{Hao Pan}
\email{haopan79@yahoo.com.cn}
\address{
Department of Mathematics, Shanghai Jiaotong University, Shanghai
200240, People's Republic of China} \thanks{This work was supported by the
National Natural Science Foundation of China (Grant No. 10771135).}

\begin{abstract}
Suppose that all primes are colored with $k$ colors. Then there
exist monochromatic primes $p_1,p_2,p_3$ such that $p_1+p_2=p_3+1$.
\end{abstract}
\maketitle

\section{Introduction}
\setcounter{equation}{0} \setcounter{Thm}{0} \setcounter{Lem}{0}
\setcounter{Cor}{0}

In \cite{GreenTao}, Green and Tao proved a celebrated result that
the primes contain arbitrarily long non-trivial arithmetic
progressions. In fact, they proved a Szemer\'edi-type
\cite{Szemeredi75} result for primes:

\medskip {\it If $A$ is a set of primes with positive relative upper
density, then $A$ contains arbitrarily long arithmetic
progressions.}

\medskip
\noindent Thus if all primes are colored with $k$ colors, then there
exist arbitrarily long monochromatic arithmetic progressions. This
is a van der Waerden-type \cite{vanderWaerden27} theorem for primes.
(The well-known van der Waerden theorem states that for any
$k$-coloring of all positive integers, there exist arbitrarily long
monochromatic arithmetic progressions.)

On the other hand, Schur's theorem \cite{Schur16} is another famous
result in the Ramsey theory for integers. Schur's theorem asserts
that for any $k$-coloring of all positive integers, there exist
monochromatic $x,y,z$ such that $x+y=z$. In this paper, we shall
prove a Schur-type theorem for primes.
\begin{Thm}
\label{t1} Suppose that all primes are arbitrarily colored with $k$
colors. Then there exist monochromatic primes $p_1,p_2,p_3$ such
that $p_1+p_2=p_3+1$.
\end{Thm}

Furthermore, motivated by the Green-Tao theorem and Theorem
\ref{t1}, we propose the following conjecture:
\begin{Conj}
Suppose that all primes are colored with $k$ colors. Then for
arbitrary $l\geq 3$, there exist monochromatic primes
$p_0,p_1,p_2,\ldots,p_l$ such that $p_1,\ldots,p_l$ form an
arithmetic progression with the difference $p_0-1$.
\end{Conj}

Theorem \ref{t1} will be proved in the next section. And our proof
uses a variant of Green's method \cite{Green05} in his proof of
Roth's theorem for primes.

\section{Proof of Theorem \ref{t1}}
\setcounter{equation}{0} \setcounter{Thm}{0} \setcounter{Lem}{0}
\setcounter{Cor}{0}

\begin{Lem} \label{schur} Suppose that the set $\{1,2,\ldots,n\}$
is split into $A_1\cup A_2\cup\cdots\cup A_k$. Then there exists a
constant $C_1(k)>0$ such that
$$
\sum_{1\leq i\leq k}|\{(x,y,z):\, x,y,z\in A_i, x+y=z\}|\geq
C_1(k)n^2
$$
if $n$ is sufficiently large.
\end{Lem}

This result is not new. In fact, Robertson and Zeilberger
\cite{RobertsonZeilberger98}, Schoen \cite{Schoen99} had showed
that if the integers from 1 to $n$ are colored with two colors,
then there exist at least $(1/22-\epsilon)n^2$ monochromatic Schur
triples $\{x,y,x+y\}$. Furthermore, Robertson and Zeilberger
\cite{RobertsonZeilberger98} also claimed that for any
$k$-coloring of $\{1,\ldots,n\}$, the number of monochromatic
Schur triples is greater than
$$
\big(\frac{1}{2^{2k-3}11}-\epsilon\big)n^2.
$$

However, for the sake of completeness, here we give a proof of
Lemma \ref{schur}. Suppose that $1,2,\ldots,n$ are colored with
$k$ colors. Let $G$ be a complete graph with the vertex set
$V=\{v_0,v_1,\ldots,v_n\}$. Then we $k$-color all edges of $G$ by
giving the edge $v_sv_t$ the color of $t-s$ for every $0\leq
s<t\leq n$. Clearly for $0\leq r<s<t\leq n$, three vertices
$v_r,v_s,v_t$ form a monochromatic triangle if and only if
$\{s-r,t-s,t-r\}$ is a monochromatic Schur triple. And it is easy
to see that one monochromatic Schur triple is corresponding to at
most $n$ monochromatic triangles. Hence Lemma \ref{schur}
immediately follows from the next lemma:
\begin{Lem} Let $G$ be a complete graph with $n$ vertices. If
all edges of $G$ are colored with $k$ colors, then there exist at
least $C_1'(k)n^3$ monochromatic triangles provided that $n$ is
sufficiently large, where $C_1'(k)>0$ is a constant only depending
on $k$.
\end{Lem}
\begin{proof}
Since $G$ is a complete graph, $G$ contains $\binom{n}{3}$
triangles. We use induction on $k$. There is nothing to do when
$k=1$. Assume that $k\geq 2$ and our assertion holds for any
smaller value of $k$. Suppose that the vertex set $V$ of $G$ is
$\{v_1,\ldots,v_n\}$. Then for every $1\leq s\leq n$, by the
pigeonhole principle, there exist vertices
$v_{t_{s,1}},\ldots,v_{t_{s,\ceil{n/k}}}$ and $1\leq c_s\leq k$
such that the edge $v_sv_{t_{s,1}},\ldots,v_sv_{t_{s,\ceil{n/k}}}$
are colored with the $c_s$-th color, where $\ceil{x}$ denotes the
smallest integer not less than $x$. Let us consider the
$\binom{\ceil{n/k}}{2}$ edges between
$v_{t_{s,1}},\ldots,v_{t_{s,\ceil{n/k}}}$. Suppose that at most
$(C_1'(k-1)/2k^3)n^2$ of these edges are colored with the $c_s$-th
color. Then by the induction hypothesis on $k-1$, the remainder
edges form at least
$$
C_1'(k-1)(n/k)^3-\frac{C_1'(k-1)}{2k^3}n^3=\frac{C_1'(k-1)}{2k^3}n^3
$$
monochromatic triangles, since one edge belongs to at most $n$
triangles.

Then we may assume that for each $1\leq s\leq n$, there exist at
least $(C_1'(k-1)/2k^3)n^2$ edges between
$v_{t_{s,1}},\ldots,v_{t_{s,\ceil{n/k}}}$ are colored with the
$c_s$-th color. Thus we get at least $(C_1'(k-1)/2k^3)n^2$
monochromatic triangles containing the vertex $v_s$. And there are
totally at least
$$
\frac{C_1'(k-1)}{6k^3}n^3
$$
monochromatic triangles, by noting that every triangles are counted
three times.
\end{proof}

\begin{Cor}
Let $A$ be a subset of $\{1,2,\ldots,n\}$ with $|A|\geq
(1-C_1(k)/6)n$. Suppose that $A$ is split into $A_1\cup
A_2\cup\cdots\cup A_k$. Then
$$
\sum_{1\leq i\leq k}|\{(x,y,z):\, x,y,z\in A_i, x+y=z\}|\geq
\frac{C_1(k)}{2}n^2
$$
provided that $n$ is sufficiently large.
\end{Cor}
\begin{proof}
Let $\bar{A}=\{1,\ldots,n\}\setminus A$. Then
\begin{align*}
&|\{(x,y,z):\, x,y,z\in A_1\cup\bar{A}, x+y=z\}|\\
\leq&|\{(x,y,z):\, x,y,z\in A_1, x+y=z\}|\\
&+|\{(x,y,z):\, \text{one of
}x,y,z\text{ lies in }\bar{A}, x+y=z\}|\\
\leq&|\{(x,y,z):\, x,y,z\in A_1, x+y=z\}|+3|\bar{A}|n.
\end{align*}
Hence by Lemma 2.1 we have
\begin{align*}
&\sum_{1\leq i\leq k}|\{(x,y,z):\, x,y,z\in A_i, x+y=z\}|\\
\geq&|\{(x,y,z):\, x,y,z\in A_1\cup\bar{A},
x+y=z\}|-3|\bar{A}|n\\
+&\sum_{2\leq i\leq k}|\{(x,y,z):\, x,y,z\in A_i,
x+y=z\}|\\
\geq&\frac{C_1(k)}{2}n^2.
\end{align*}

\end{proof}

Let $\P$ denote the set of all primes. Assume that $\P=P_1\cup
P_2\cup\cdots\cup P_k$, where $P_i\cap P_j=\emptyset$ for $1\leq
i<j\leq k$. Let $w=w(n)$ be a function tending sufficiently slowly
to infinity with $n$ (e.g., we may choose
$w(n)=\floor{\frac{1}{4}\log\log n}$), and let
$$
W=\prod_{\substack{p\in\P\\ p\leqslant w(n)}}p.
$$
Clearly  we have $W\leqslant\log n$ for sufficiently large $n$.
Let
$$
\kappa=\frac{C_1(k)}{10000k}.
$$
In view of the well-known Siegel-Walfisz theorem, we may assume that
$n$ is sufficiently large so that
$$
\sum_{\substack{x\in\P\cap[1,n]\\ x\equiv 1\pmod{W}}}\log
x\geq(1-\kappa)\frac{n}{\phi(W)},
$$
where $\phi$ is the Euler totient function. Let $M=n/W$ and $N$ be a
prime in the interval $[(2+\kappa)M,(2+2\kappa)M]$. (Thanks to the
prime number theorem, such prime $N$ always exists whenever $M$ is
sufficiently large.) Define
$$
\lambda_{b,W,N}(x)=\begin{cases} \phi(W)\log(Wx+b)/WN&\text{ if
}x\leq N\text{ and }Wx+b\text{ is prime},\\
0&\text{otherwise}.
\end{cases}
$$
Let
$$
A_0=\{1\leq x\leq M:\, Wx+1\in\P\}
$$
and
$$
A_i=\{1\leq x\leq M:\, Wx+1\in P_i\}
$$
for $1\leq i\leq k$. Define $$a_i(x)=\1_{A_i}(x)\lambda_{1,W,N}(x)$$
for $0\leq i\leq k$, where we set $\1_A(x)=1$ if $x\in A$ and $0$
otherwise. Clearly we have $a_0=a_1+\cdots+a_k$ and
$$
\sum_{x}a_0(x)=\sum_{1\leq x\leq M}\lambda_{1,W,N}(x)\geq
(1-\kappa)\frac{M}{N}\geq \frac{1}{2}-3\kappa.
$$

Below we consider $A_0,A_1,\ldots,A_k$ as the subsets of
$\Z_N=\Z/N\Z$. Since $M<N/2$, if there exist $x,y,z\in A_i$ such
that $x+y=z$ in $\Z_N$, then we have $p_1+p_2=p_3+1$ in $\Z$, where
$p_1=Wx+1\in P_i, p_2=Wy+1\in P_i, p_3=Wz+1\in P_i$. For a
complex-valued function $f$ over $\Z_N$, define $\tilde{f}$ by
$$
\tilde{f}(r)=\sum_{x\in\Z_N}f(x)e(-xr/N),
$$
where $e(x)=e^{2\pi\sqrt{-1}x}$. And for two functions $f, g$,
define
$$
(f*g)(x)=\sum_{y\in\Z_N}f(y)g(x-y).
$$
It is easy to check that $(f*g)\,\tilde{}=\tilde{f}\tilde{g}$. Let
$0<\delta, \epsilon<1/2$ be two sufficiently small real numbers
which will be chosen later. Let
$$
R=\{r\in\Z_N:\, \max_{1\leq i\leq
k}|\tilde{a_i}(r)|\geqslant\delta\}.
$$
and
$$
B=\{x\in\Z_N:\, x\in[-\kappa N,\kappa N],\
\|xr/N\|\leqslant2\epsilon\text{ for all }r\in R\},
$$
where $ \|x\|=\min_{z\in\Z}|x-z|$. Here our definition of $B$ is
slightly different from Green's one in \cite[Page 1629]{Green05}. As
we shall see later, this modification is the key of our proof.
\begin{Lem}
\label{bohr}
$$
|B|\geq\epsilon^{|R|}\kappa N.
$$
\end{Lem}
\begin{proof}
Assume that $R=\{r_1,r_2,\ldots,r_{m}\}$. Let $d$ be the greatest
integer not exceeding $1/\epsilon$. Clearly we have $1/d\leq
2\epsilon$ since $\epsilon<1/2$. Let
$$
G_{t_1,\ldots,t_m}=\{-\kappa N/2\leq x\leq\kappa N/2:\,
t_j/d\leq\{xr_j/N\}<(t_j+1)/d\text{ for }1\leq j\leq m\},
$$
where $\{\alpha\}$ denotes the fractional part of $\alpha$.  Clearly
$$
\sum_{0\leq t_1,\ldots,t_m\leq d-1}|G_{t_1,\ldots,t_m}|=\kappa N.
$$
Hence there exists a term of $(t_1,\ldots,t_m)$ such that
$$
|G_{t_1,\ldots,t_m}|\geq d^{-m}\kappa N\geq \epsilon^{m}\kappa N.
$$
For any given $x_0\in G_{t_1,\ldots,t_{m}}$, when $x\in
G_{t_1,\ldots,t_{m}}$, we have $x-x_0\in [-\kappa N,\kappa N]$ and
$$
\|(x-x_0)r_j/N\|\leq|\{xr_j/N\}-\{x_0r_j/N\}|\leq1/d\leq2\epsilon
$$
for $1\leq j\leq m$. So $G_{t_1,\ldots,t_{m}}\subseteq x_0+B$. This
completes the proof.
\end{proof}
\begin{Lem}
\label{lambda}
$$
\sup_{r\not=0}|\tilde{\lambda}_{b,W,N}(r)|\leq 2\log\log w/w
$$
provided that $w$ is sufficiently large.
\end{Lem}
\begin{proof}
This is Lemma 6.2 of \cite{Green05}.
\end{proof}
Let $\beta=\1_{B}/|B|$ and $a_i'=a_i*\beta*\beta$ for $0\leq i\leq
k$.
\begin{Lem}
\label{upper} Suppose that $\epsilon^{|R|}\geq \kappa^{-2}\log\log
w/w$. Then we have
$$
\sup_{x\in\Z_N}a_0'(x)\leq\frac{1+3\kappa}{N}.
$$
\end{Lem}
\begin{proof} We have
\begin{align*}
a_0'(x)=&a_0*\beta*\beta(x)\\
\leq&\lambda_{1,W,N}*\beta*\beta(x)\\
=&N^{-1}\sum_{r\in\Z_N}\tilde{\lambda}_{1,W,N}(r)\tilde{\beta}(r)^2e(xr/N)\\
\leq&N^{-1}\tilde{\lambda}_{1,W,N}(0)\tilde{\beta}(0)^2+N^{-1}\sup_{r\not=0}|\tilde{\lambda}_{1,W,N}(r)|\sum_{r\in\Z_N}|\tilde{\beta}(r)|^2\\
=&N^{-1}\tilde{\lambda}_{1,W,N}(0)+\sup_{r\not=0}|\tilde{\lambda}_{1,W,N}(r)|\sum_{r\in\Z_N}|\beta(r)|^2\\
\leq&\frac{1+\kappa}{N}+\frac{2\log\log w}{w|B|},
\end{align*}
where Lemma \ref{lambda} is applied in the last step. Thus Lemma
\ref{upper} immediately follows from Lemma \ref{bohr}.
\end{proof}

\begin{Lem}[Bourgain \cite{Bourgain89,Bourgain93}, Green \cite{Green05}]
\label{bourgain} Let $\rho>2$. For any function $f:\,\Z_N\to\C$,
$$
\sum_{r\in\Z_N}|(f\lambda_{b,W,N})\,\tilde{}(r)|^\rho\leq
C_2(\rho)\bigg(\sum_{x=1}^N|f(x)|^2\lambda_{b,W,N}(x)\bigg)^{\frac{\rho}{2}}
$$
where $C_2(\rho)$ is a constant only depending on $\rho$.
\end{Lem}
\begin{proof}
This is an immediate consequence of Theorem 2.1 and Lemma 6.5 of
\cite{Green05}.
\end{proof}
By Lemma \ref{bourgain}, we have
$$
\sum_{r\in\Z_N}|\tilde{a}_i(r)|^\rho\leq C_2(\rho)
$$
for $\rho>2$ and $1\leq i\leq k$. In particular,
$$
\sum_{r\in R}\delta^3\leq
\sum_{r\in\Z_N}\bigg(\sum_{i=1}^k|\tilde{a}_i(r)|^3\bigg)\leq
C_2(3)k,
$$
which implies that $|R|\leq C_2(3)\delta^{-3}k$.

\begin{Lem}
\label{beta} For each $r\in R$,
$$
|1-\tilde{\beta}(r)^4\tilde{\beta}(-r)^2|\leq 384\epsilon^2.
$$
\end{Lem}
\begin{proof}
By the definition of $B$, we have
\begin{align*}
|1-\tilde{\beta}(r)|=\frac{1}{|B|}\bigg|\sum_{x\in
B}(1-e(-xr/N))\bigg| \leq4\pi\sup_{x\in B}\|xr/N\|^2\leq
64\epsilon^2.
\end{align*}
So
\begin{align*}
|1-\tilde{\beta}(r)^4\tilde{\beta}(-r)^2|=&\bigg|\sum_{j=0}^3\tilde{\beta}(r)^j(1-\tilde{\beta}(r))+\tilde{\beta}(r)^4\sum_{j=0}^1\tilde{\beta}(-r)^j(1-\tilde{\beta}(-r))\bigg|\\
\leq&384\epsilon^2.
\end{align*}
by noting that $|\tilde{\beta}(r)|\leq\tilde{\beta}(0)=1$.
\end{proof}
\begin{Lem}
\label{difference} For $1\leq i\leq k$,
$$
\bigg|\sum_{i=1}^k\sum_{\substack{
x,y,z\in\Z_N\\x+y=z}}a_i(x)a_i(y)a_i(z)-\sum_{i=1}^k\sum_{\substack{
x,y,z\in\Z_N\\x+y=z}}a_i'(x)a_i'(y)a_i'(z)\bigg|\leq
\frac{C_3k^2}{N}(\epsilon^2\delta^{-3}+\delta^{\frac{1}{3}}),
$$
where $C_3$ is an absolute constant.
\end{Lem}
\begin{proof}
Clearly
$$
\sum_{\substack{
x,y,z\in\Z_N\\x+y=z}}f_1(x)f_2(y)f_3(z)=N^{-1}\sum_{r\in\Z_N}\tilde{f}_1(r)\tilde{f}_2(r)\tilde{f}_3(-r).
$$
Hence
\begin{align*}
&\sum_{i=1}^k\sum_{\substack{
x,y,z\in\Z_N\\x+y=z}}a_i(x)a_i(y)a_i(z)-\sum_{i=1}^k\sum_{\substack{
x,y,z\in\Z_N\\x+y=z}}a_i'(x)a_i'(y)a_i'(z)\\
=&N^{-1}\sum_{i=1}^k\sum_{r\in\Z_N}\tilde{a}_i(r)^2\tilde{a}_i(-r)(1-\tilde{\beta}(r)^4\tilde{\beta}(-r)^2).
\end{align*}
By Lemma \ref{beta},
\begin{align*}
&\bigg|\sum_{i=1}^k\sum_{r\in
R}\tilde{a}_i(r)^2\tilde{a}_i(-r)(1-\tilde{\beta}(r)^4\tilde{\beta}(-r)^2)\bigg|\\
\leq&384\epsilon^2k|R|\sup_{r}\max_{1\leq i\leq
k}|\tilde{a}_i(r)|^3\\\leq& 384C_2(3)\epsilon^2\delta^{-3}k^2,
\end{align*}
since $|\tilde{a}_i(r)|\leq\tilde{a}_i(0)\leq 1$. On the other hand,
by the H\"older inequality, we have
\begin{align*}
&\bigg|\sum_{i=1}^k\sum_{r\not\in
R}\tilde{a}_i(r)^2\tilde{a}_i(-r)(1-\tilde{\beta}(r)^4\tilde{\beta}(-r)^2)\bigg|\\
\leq&2\sum_{i=1}^k\sum_{r\not\in
R}|\tilde{a}_i(r)|^2|\tilde{a}_i(-r)|\\
\leq&2\sup_{r\not\in R}\max_{1\leq i\leq
k}|\tilde{a}_i(r)|^{\frac{1}{3}}\bigg(\sum_{i=1}^k\sum_{r}|\tilde{a}_i(r)|^{\frac{5}{2}}\bigg)^{\frac{2}{3}}
\bigg(\sum_{i=1}^k\sum_{r}|\tilde{a}_i(-r)|^{3}\bigg)^{\frac{1}{3}}\\
\leq&2C_2(5/2)^{\frac{2}{3}}C_2(3)^{\frac{1}{3}}\delta^{\frac13}k.
\end{align*}
We choose
$C_3=384C_2(3)+2C_2(5/2)^{\frac{2}{3}}C_2(3)^{\frac{1}{3}}$, then
the Lemma follows.
\end{proof}

Define
$$
X=\{x\in\Z_N:\, a_0'(x)\geq \frac{\kappa}{N} \}.
$$
Then by Lemma \ref{upper}, we have
$$
\frac{1+3\kappa}{N}|X|+\frac{\kappa}{N}(N-|X|)\geq
\sum_{x\in\Z_N}a_0'(x)=\sum_{x\in\Z_N}a_0(x)\geq\frac{1}{2}-3\kappa.
$$
It follows that
$$
|X|\geq\big(\frac{1}{2}-6\kappa\big)N.
$$
Notice that $\supp(a_i)\subseteq[1,M]$ and
$\supp(\beta)\subseteq[-\kappa N,\kappa N]$, where $$
\supp(f)=\{x\in\Z_N:\, f(x)\not=0\}.
$$
Hence $$\supp(a_i')=\supp(a_i*\beta*\beta)\subseteq[-2\kappa
N,M+2\kappa N]$$ for $0\leq i\leq k$. Thus we have
$$
X\subseteq\supp(a_0')\subseteq[-2\kappa N,M+2\kappa N].
$$ Let $A_0'=X\cap[1,M]$. Then
$$
|A_0'|\geq|X|-4\kappa N\geq(1-20\kappa)M,
$$
by recalling that $(2+\kappa)M\leq N\leq (2+2\kappa)M$. Since
$$
a_0'=a_0*\beta*\beta=(a_1+\cdots+a_k)*\beta*\beta=a_1'+\cdots+a_k',
$$
we have
$$
\max_{1\leq i\leq k}a_i'(x)\geq\frac{\kappa}{kN}
$$
for each $x\in A_0'$. Let
$$
X_i=\{x\in A_0':\, a_i'(x)=\max_{1\leq i\leq k}a_i'(x)\}.
$$
Clearly $A_0'=X_1\cup\cdots\cup X_k$. Let $A_1'=X_1$ and
$$
A_i'=X_i\setminus\bigg(\bigcup_{j=1}^{i-1}X_j\bigg)
$$
for $2\leq i\leq k$. Then $A_1',\ldots,A_k'$ form a partition of
$A_0'$. Furthermore, for $1\leq i\leq k$ and each $x\in A_i'$, we
have
$$
a_i'(x)\geq\frac{\kappa}{kN}.
$$
Thus by Corollary \ref{schur} and Lemma \ref{difference}
\begin{align*}
\sum_{i=1}^k\sum_{\substack{ x,y,z\in\Z_N\\x+y=z}}a_i(x)a_i(y)a_i(z)
\geq&\sum_{i=1}^k\sum_{\substack{
x,y,z\in\Z_N\\x+y=z}}a_i'(x)a_i'(y)a_i'(z)-\frac{C_3k^2}{N}(\epsilon^2\delta^{-3}+\delta^{\frac{1}{3}})\\
\geq&\sum_{i=1}^k\sum_{\substack{
x,y,z\in A_i'\\x+y=z}}\bigg(\frac{\kappa}{kN}\bigg)^3-\frac{C_3k^2}{N}(\epsilon^2\delta^{-3}+\delta^{\frac{1}{3}})\\
\geq&\bigg(\frac{\kappa}{kN}\bigg)^3\frac{C_1(k)M^{2}}{2}-\frac{C_3k^2}{N}(\epsilon^2\delta^{-3}+\delta^{\frac{1}{3}})\\
\end{align*}
Finally, we may choose sufficiently small $\delta$ and $\epsilon$
with
$$
\epsilon^{-C_2(3)\delta^{-3}k}\geq \kappa^{-2}\log\log w/w
$$
such that
$$
\epsilon^2\delta^{-3}+\delta^{\frac{1}{3}}\leq\frac{C_1(k)\kappa^3}{24C_3k^5},
$$
whenever $N$ is sufficiently large. Thus
\begin{align*}
\sum_{i=1}^k\sum_{\substack{ x,y,z\in\Z_N\\x+y=z}}a_i(x)a_i(y)a_i(z)
\geq\frac{C_1(k)\kappa^3M^{2}}{2k^3N^3}-\frac{C_1(k)\kappa^3}{24k^3N}\geq
\frac{C_1(k)\kappa^3}{12k^3N}-\frac{C_1(k)\kappa^3}{24k^3N}>0.
\end{align*}
This completes the proof. \qed

\begin{Rem}
Notice that
$$
\sum_{i=1}^k\sum_{\substack{ x,z\in\Z_N\\
2x=z}}a_i(x)^2a_i(z)=O\bigg(\frac{k\phi(W)^3\log(WN+1)^3}{W^3N^2}\bigg)=o(N^{-1}).
$$
Hence in fact there exist three distinct monochromatic primes
$p_1,p_2,p_3$ satisfying $p_1+p_2=p_3+1$.
\end{Rem}

\end{document}